\documentclass[11pt]{article}

\usepackage[a4paper,top=30mm,bottom=30mm,left=31mm,right=31mm]{geometry}
\usepackage[T1]{fontenc}
\usepackage{amsmath,amsthm}
\usepackage{newtxtext,newtxmath}
\usepackage{microtype}
\usepackage{mathtools}
\usepackage{enumitem}
\usepackage{booktabs}
\usepackage{xcolor}
\usepackage{titlesec}
\usepackage{fancyhdr}
\usepackage[hidelinks]{hyperref}
\usepackage[nameinlink,noabbrev]{cleveref}

\allowdisplaybreaks[2]
\setlist{nosep,leftmargin=2.2em}
\titleformat{\section}{\large\bfseries}{\thesection.}{0.55em}{}
\titleformat{\subsection}{\normalsize\bfseries}{\thesubsection.}{0.55em}{}
\titlespacing*{\section}{0pt}{2.15ex plus .6ex minus .2ex}{1.05ex}
\titlespacing*{\subsection}{0pt}{1.65ex plus .4ex minus .2ex}{.7ex}

\newtheorem{theorem}{Theorem}[section]
\newtheorem{proposition}[theorem]{Proposition}
\newtheorem{lemma}[theorem]{Lemma}
\newtheorem{corollary}[theorem]{Corollary}

\theoremstyle{definition}

\theoremstyle{remark}
\newtheorem{remark}[theorem]{Remark}

\newcommand{\K}{\mathbb K}
\newcommand{\N}{\mathbb N}
\newcommand{\whot}{\widehat{\otimes}}
\newcommand{\tp}{\mathbin{\whot_{\!\pi}}}
\newcommand{\tpr}[1]{\mathbin{\whot_{\!\pi_{#1}}}}
\newcommand{\ce}{\stackrel{c}{\hookrightarrow}}

\newcommand{\norm}[1]{\left\lVert #1\right\rVert}
\newcommand{\abs}[1]{\left|#1\right|}
\newcommand{\doi}[1]{\href{https://doi.org/#1}{doi:#1}}

\begin{document}

\begin{center}
{\LARGE\bfseries Isomorphic nonembeddings of Lorentz spaces under tensor products}\par
\vspace{1.1em}
{\large Ying Xu$^{1}$ Zhiyao Fang$^{1}$ and Qi Liu$^{1*}$}\par
\vspace{0.8em}
$^{1*}$School of Mathematics and Statistics, Anqing Normal University, Anqing, 246133,\par
Anhui, China.\par
\vspace{1.0em}
\textit{Corresponding author(s).} E-mail(s): \href{mailto:liuq67@aqnu.edu.cn}{liuq67@aqnu.edu.cn};\par
\textit{Contributing authors:} \href{mailto:Y26220046@stu.aqnu.edu.cn 
	 fzy73@stu.aqnu.edu.cn}{Y26220046@stu.aqnu.edu.cn \quad fzy73@stu.aqnu.edu.cn};\par
\end{center}

\begin{abstract}
	The Lorentz quasi-norm on $L_{p,q}(0,1)$ is
	\[
	\|f\|_{p,q}
	=\left(\int_0^1 [t^{1/p}f^*(t)]^q\,\frac{dt}{t}\right)^{1/q},
	\]
	where $f^*$ is the decreasing rearrangement of $|f|$.
	Motivated by the isomorphic classification of Lorentz spaces
	\cite{KuryakovSukochev2015,HuangSukochev2021}, we ask whether
	$\ell_{p,q}\not\hookrightarrow L_{p,q}(0,1)$ persists under projective
	tensoring with a common factor. For $1<p<\infty$, $1\le q<\infty$,
	and $p\ne q$, we obtain the following results.
	
	\begin{enumerate}[label=\textup{(\alph*)}, leftmargin=*]
		\item Nonzero finite-dimensional factors preserve the nonembedding,
		as do nonzero $\ell_1$-FDD factors when $q>1$.
		\item A single factor $U$ makes both tensor products isomorphic to $U$
		throughout this range. For each pair with $q>1$, there are two
		nonisomorphic separable absorbers with the metric approximation
		property.
		\item For $1<q<p<2$, a scalar $L_1(\Omega,\mu)$ factor over a standard
		probability space yields an embedding exactly when the nonatomic
		part of $\mu$ has positive measure.
	\end{enumerate}
	
	For $1\le p<\infty$ and $0<q<1$, we construct common absorbers under
	$r$-projective tensoring after equivalent renorming.
\end{abstract}



\section{Introduction}

Lorentz spaces are among the most economical models in which measure theory, interpolation and Banach-space geometry meet.  They retain the distributional scale of $L_p$ while separating the integrability index from the summability index, and this extra parameter records geometric information that is invisible in the Lebesgue scale; see Lorentz's original paper \cite{Lorentz1950}, Hunt's systematic treatment \cite{Hunt1966}, and Bennett--Sharpley \cite{BennettSharpley1988}.  For $0<p<\infty$ and $0<q<\infty$, the function-space model $L_{p,q}(0,1)$ is equipped with the Lorentz quasi-norm
\begin{equation}\label{eq:lorentz-quasinorm}
    \|f\|_{p,q}
    =\left(\int_0^1 [t^{1/p}f^*(t)]^q\,\frac{dt}{t}\right)^{1/q},
\end{equation}
where $f^*$ is the decreasing rearrangement of $|f|$.  Thus, $p$ is the primary integrability parameter, while $q$ is the secondary rearrangement summability parameter.  For fixed $p$ and $q$, the sequence model $\ell_{p,q}$ and the function model $L_{p,q}(0,1)$ live on opposite sides of a measure-theoretic divide.  The former is purely atomic, whereas the latter is nonatomic.  Although their rearrangement formulas have the same formal shape, their linear structures need not agree.

The sharp manifestation relevant here is the isomorphic nonembedding between the atomic and nonatomic Lorentz models. Here and throughout, the notation $X\hookrightarrow Y$ means that $X$ is
linearly isomorphic to a subspace of $Y$; equivalently, there exists a bounded
linear embedding $T:X\to Y$ whose inverse on $T(X)$ is bounded. We write
$X\not\hookrightarrow Y$ when no such isomorphic embedding exists. The notation
$X\cong Y$ will be used for a Banach-space isomorphism between $X$ and $Y$,
and does not imply that the two spaces are isometric. The fundamental obstruction considered in this work is
\begin{equation}\label{eq:base-nonembedding}
    \ell_{p,q}\not\hookrightarrow L_{p,q}(0,1),
\end{equation}
throughout
\begin{equation}\label{eq:banach-range}
    1<p<\infty,\qquad 1\le q<\infty,\qquad p\ne q.
\end{equation}
This relation expresses the failure of the atomic Lorentz sequence space to embed isomorphically into its corresponding nonatomic function-space model. We take \eqref{eq:base-nonembedding} as the underlying geometric obstruction and ask whether it remains stable or can be eliminated after tensor stabilization.

The result follows from the work of Kuryakov--Sukochev \cite{KuryakovSukochev2015}, together with the quadratic endpoint analysis of Sadovskaya--Sukochev \cite{SadovskayaSukochev2018}. It belongs to the isomorphic classification program for Lorentz spaces on resonant measure spaces, later completed in greater generality by Huang--Sukochev \cite{HuangSukochev2021}; a noncommutative counterpart has recently been obtained by Huang, Sadovskaya, Sukochev and Zanin \cite{HuangEtAl2025}.

For Banach spaces $E$ and $X$, let $E\tp X$ denote the completed projective tensor product.  The same-factor problem is
\begin{equation}\label{eq:same-factor}
  \ell_{p,q}\tp X \ \hookrightarrow\ \ L_{p,q}(0,1)\tp X.
\end{equation}
This is not a cancellation problem in a naive algebraic sense.  Projective tensoring can preserve a rigid subspace obstruction, but it can also enlarge both spaces so much that the original distinction disappears.  The right invariant is therefore not the dimension of $X$ but the geometry that $X$ forces into the tensor product.

The paper develops this principle in two categories. In the Banach category, we use the classical projective norm, complemented universality, and Pełczyński's decomposition method, while beyond local convexity we employ the topological tensor products associated with quasi-Banach tensor constructions. These approaches allow us to analyze how decomposition structures, approximation properties, and dual information determine whether Lorentz nonembedding phenomena are preserved or eliminated under tensor products.

The projective tensor product viewpoint places the present problem within the broader framework of Banach space geometry. Tensor norms and operator ideals provide a systematic language for describing the interaction between geometric properties of individual spaces after tensorization. A fundamental difficulty is that tensor products do not have a uniform effect on geometric distinctions: some tensor factors preserve the atomic--nonatomic separation, whereas others absorb the difference completely. Therefore, the central issue is to identify the stability mechanism governing Lorentz nonembedding under tensor operations.

To analyze this phenomenon, we combine classical tools including the bounded approximation property (BAP), complemented universality, and Pełczyński's decomposition method \cite{Pelczynski1960,Pelczynski1969,Pelczynski1971}. These techniques distinguish factors that preserve the original obstruction from absorbing factors that eliminate it through tensor identities. We further extend the analysis to the quasi-Banach setting, where the lack of local convexity introduces additional difficulties. By exploiting quasi-Banach tensor constructions and the role of nonzero continuous functionals, we identify the structural conditions required for absorption arguments. 

The preceding results describe the Lorentz nonembedding phenomenon at the
level of individual spaces. A natural question is whether this obstruction
remains stable under tensor amplification. In particular, one may ask whether
the projective tensor product preserves the atomic--nonatomic Lorentz
distinction, or whether a suitable tensor factor can fundamentally change this
geometry.

The main contribution of this paper is to establish a tensor-stability framework 
for the atomic--nonatomic Lorentz nonembedding phenomenon. We show that the 
obstruction in \eqref{eq:base-nonembedding} is determined not only by the original 
Lorentz pair, but also by the geometric structure imposed by the tensor factor: 
{finite-dimensional factors preserve the nonembedding, as do $\ell_1$-FDD factors in the reflexive range,} 
whereas universal BAP absorbers can eliminate this obstruction. Thus, the paper 
identifies the mechanism by which projective tensor products either retain or 
remove the Lorentz geometric distinction.

Throughout, an \emph{isomorphic embedding} means a bounded linear injection bounded below, and $E\ce F$ means that $E$ is isomorphic to a complemented subspace of $F$.  The paper is organized as follows. Section~2 collects the preliminary tensor, approximation, and localization tools used throughout. Section~3 develops the Banach stabilization and rigidity theory, {including common-factor absorption, nonuniqueness among separable MAP factors, and the finite-dimensional and $\ell_1$-FDD stability results.} Section~4 treats the quasi-Banach extension and the trivial-dual obstruction. Section~5 records further consequences and closes with the resulting stability--absorption picture.

{A further main result is Theorem~\ref{thm:map-nonunique-factors}: for every fixed pair with $1<p,q<\infty$, $p\ne q$, two nonisomorphic separable MAP factors absorb both Lorentz spaces. This differs from the orbit uniqueness in Theorem~\ref{thm:bap-rigidity}. These constructions and the preservation results describe specified classes of factors; no exhaustive classification of all common factors is asserted.}

\section{Preliminary tools}\label{sec:prelim}

We collect here the structural tools used in the Banach and quasi-Banach arguments.  The purpose of this section is preparatory: the main stabilization theorems are stated only after the tensor identities, complemented embeddings, and localization lemmas on which their proofs depend.

The projective tensor product is the natural framework for our questions because
the arguments below rely on complemented decompositions, factorization
properties, and universality phenomena, all of which interact directly with
the projective norm.  Thus, throughout the paper, tensor stability refers to
the behavior of the Lorentz nonembedding obstruction under projective tensor
amplification.

We work over $\K=\mathbb R$ or $\mathbb C$.  For Banach spaces $E$ and $F$, the projective norm on $E\otimes F$ is
\[
 \norm{u}_\pi=\inf\left\{\sum_{k=1}^n\norm{x_k}\norm{y_k}:
 u=\sum_{k=1}^n x_k\otimes y_k\right\},
\]

and $E\tp F$ denotes its completion.  We use standard projective-tensor facts from \cite{DefantFloret1993,Ryan2002,DiestelFourieSwart2008}: functoriality; symmetry and associativity; and the canonical isometry
\begin{equation}\label{eq:l1-distrib}
 E\tp\left(\bigoplus_nF_n\right)_{\ell_1}
 \cong
 \left(\bigoplus_n E\tp F_n\right)_{\ell_1}.
\end{equation}
We also use without reproving the classical Pe\l czy\'nski decomposition principle \cite{Pelczynski1960}: if $X\ce Y$, $Y\ce X$ and $X\cong\ell_1(X)$, then $X\cong Y$.  The BAP passes to complemented subspaces and is stable under the projective tensor product of two BAP spaces; see \cite{DefantFloret1993,Ryan2002}.  Every space with a Schauder basis has the BAP.

\begin{lemma}
\label{lem:rank-one-tensor-retraction}
Let $A$ and $X$ be nonzero Banach spaces. Then $A$ is isometric to a
$1$-complemented subspace of $A\tp X$.
\end{lemma}

\begin{proof}
Choose $x_0\in X$ and $x^*\in X^*$ such that $x^*(x_0)=1$.  The map
$a\mapsto a\otimes x_0$ gives an isometric embedding of $A$ into
$A\tp X$, while the operator
\[
a\otimes x\longmapsto x^*(x)a
\]
extends to a contractive projection onto this copy of $A$.
\end{proof}

The absorption and rigidity phenomena in the Banach category are closely connected with the existence of universal complemented structures. To describe this mechanism, we recall the classical construction of Pe\l czy\'nski. He constructed a separable space $P$ with a Schauder basis that is complementably universal for separable spaces with a basis \cite{Pelczynski1969}.  Combined with his theorem that every separable BAP space is complemented in a space with a basis \cite{Pelczynski1971}, one obtains a complementably universal model for the entire separable BAP class. Set
\begin{equation}\label{eq:U-def}
  U=\left(\bigoplus_{n=1}^{\infty}P\right)_{\ell_1}.
\end{equation}
Then $U$ has a Schauder basis and the BAP, is complementably universal for separable BAP spaces, and $U\cong\ell_1(U)$.

\begin{proposition}\label{prop:U-properties}
The space $U$ is separable, has a Schauder basis and the BAP, is complementably universal for separable BAP spaces, and satisfies
\[
        U\cong\ell_1(U).
\]
\end{proposition}

\begin{proof}
All assertions except self-similarity are inherited from $P$ and the $\ell_1$-sum construction.  For the last one, expand both sides as countable $\ell_1$-sums of copies of $P$ and reindex $\N\times\N$ by $\N$.
\end{proof}

\begin{lemma}\label{lem:finite-amplification}
Let $1\le m<\infty$.
\begin{enumerate}[label=\textup{(\roman*)}]
\item If $F$ is an $m$-dimensional Banach space and $E$ is Banach, then $E\tp F$ is isomorphic to a finite direct sum of $m$ copies of $E$.
\item For every $0<p,q<\infty$, the quasi-Banach spaces $(\ell_{p,q})^m$ and $\ell_{p,q}$ are isomorphic, and likewise $(L_{p,q}(0,1))^m\cong L_{p,q}(0,1)$.
\end{enumerate}
\end{lemma}

\begin{proof}
For \textup{(i)}, choose a basis $(f_j)_{j=1}^m$ of $F$ and its biorthogonal functionals.  The coordinate map $(x_j)\mapsto\sum_jx_j\otimes f_j$ and the slice maps $\mathrm{Id}_E\otimes f_j^*$ give inverse estimates.  For \textup{(ii)}, interlace $m$ coordinate copies in the sequence case; in the function case partition $(0,1)$ into $m$ equal intervals and use dilation on each band.  The Lorentz rearrangement functional scales by $m^{-1/p}$ on each band.  No local convexity is used, and all quasi-norms on a fixed finite product are equivalent.
\end{proof}

The identity \eqref{eq:l1-distrib} turns tensoring with an $\ell_1$-FDD factor into an $\ell_1$-sum of finite amplifications.  Reflexive subspaces of such sums cannot remain genuinely spread over infinitely many coordinates.  We use the following standard localization fact in an explicit varying-coordinate form.

\begin{lemma}\label{lem:uniform-tails}
Let $(B_n)_{n\ge1}$ be Banach spaces and
\[
      W=\left(\bigoplus_{n=1}^{\infty}B_n\right)_{\ell_1}.
\]
If $K\subset W$ is relatively weakly compact and $P_N$ is the projection onto the first $N$ coordinates, then
\[
   \lim_{N\to\infty}\sup_{z\in K}\norm{(I-P_N)z}=0.
\]
\end{lemma}

\begin{proof}
Suppose not.  Then, after passing to a weakly convergent sequence in $K$ and subtracting its weak limit, one obtains a weakly null sequence $(z_k)$, a number $\eta>0$, and successive disjoint coordinate intervals $I_k$ such that
\[
    \sum_{n\in I_k}\norm{(z_k)_n}>\eta
\]
while the contributions of $z_k$ on earlier and later chosen blocks are, after a standard gliding-hump selection, smaller than $\eta/8$.  For each $n\in I_k$, choose $b_n^*\in B_n^*$ of norm at most one that norms $(z_k)_n$ (with a unimodular correction over $\mathbb C$).  Since
\[
    W^*=\left(\bigoplus_{n=1}^{\infty}B_n^*\right)_{\ell_\infty},
\]
the block functionals combine into a single $\Phi\in W^*$ with $\norm{\Phi}\le1$ and $\abs{\Phi(z_k)}\ge\eta/2$ for every $k$.  This contradicts $z_k\rightharpoonup0$.
\end{proof}

\begin{corollary}\label{cor:compression}
Let $R$ be reflexive and $T:R\to W$ an isomorphic embedding, with $W$ as in Lemma \ref{lem:uniform-tails}.  Then $P_NT$ is an isomorphic embedding for some $N$.
\end{corollary}

\begin{proof}
Choose $c>0$ with $\norm{Tx}\ge c\norm{x}$.  Since $T(B_R)$ is relatively weakly compact, Lemma \ref{lem:uniform-tails} gives $N$ for which $\norm{(I-P_N)T}<c/2$.  Hence
\[
    \norm{P_NTx}\ge \frac{c}{2}\norm{x}.
\]
\end{proof}

\section{Preservation of Lorentz Nonembeddings under Projective Tensor Products}\label{sec:banach-absorption}

We now turn to the main Banach-space results.  The universal absorber describes how tensorization can absorb the original
Lorentz distinction, while the finite-dimensional and $\ell_1$-FDD results
identify factor classes for which the nonembedding obstruction remains
stable.

The introduction of a universal absorber is not merely a construction of
a special example.  It identifies a structural mechanism showing that
projective tensorization may alter the original Lorentz obstruction.  In this
sense, the following result describes the absorption side of the
{stability and absorption results.}

\begin{theorem}\label{thm:universal-absorption}
If $E\ne\{0\}$ is a separable Banach space with the bounded approximation property, then
\[
       E\tp U\cong U.
\]
Consequently, for every pair satisfying \eqref{eq:banach-range},
\[
  \ell_{p,q}\tp U\cong U\cong L_{p,q}(0,1)\tp U.
\]
\end{theorem}

\begin{proof}
The space $E\tp U$ is separable and has the BAP.  By complementable universality,
\[
       E\tp U\ce U.
\]
Because $E\ne0$, the rank-one retraction gives the reverse complemented embedding
\[
       U\ce E\tp U.
\]
Proposition \ref{prop:U-properties} gives $U\cong\ell_1(U)$, so Pe\l czy\'nski decomposition yields $E\tp U\cong U$.

For \eqref{eq:banach-range}, both Lorentz spaces are nonzero, separable and have the BAP; this follows, for example, from their standard finite-rank approximation structure in the separable Lorentz range; see \cite{BennettSharpley1988}.  Applying the first assertion to $E=\ell_{p,q}$ and $E=L_{p,q}(0,1)$ proves the stated collapse.
\end{proof}

The absorber is not merely an existence device.  Its principal tensor orbit has an exact regularity test.

\begin{theorem}\label{thm:bap-rigidity}
Let $Z\ne\{0\}$ be separable.  The following assertions are equivalent:
\begin{enumerate}[label=\textup{(\roman*)}]
\item $Z$ has the BAP;
\item $U\tp Z$ has the BAP;
\item $U\tp Z$ has a Schauder basis;
\item $U\tp Z\cong U$.
\end{enumerate}
Hence $[U]$ is the unique BAP isomorphism class among the separable members of the orbit $\{U\tp Z:Z\ne0\}$.
\end{theorem}

\begin{proof}
If $Z$ has the BAP, Theorem \ref{thm:universal-absorption}, applied to $Z$, gives $U\tp Z\cong U$ by symmetry of $\tp$.  Hence \textup{(i)}$\Rightarrow$\textup{(iv)}.  Since $U$ has a Schauder basis, \textup{(iv)}$\Rightarrow$\textup{(iii)}, and every space with a Schauder basis has the BAP, so \textup{(iii)}$\Rightarrow$\textup{(ii)}.

Assume \textup{(ii)}.  As $U\ne0$, the rank-one retraction makes $Z$ a complemented subspace of $U\tp Z$.  The BAP passes to complemented subspaces; hence $Z$ has the BAP.  Thus \textup{(ii)}$\Rightarrow$\textup{(i)} and all four assertions are equivalent.
\end{proof}

{Further tensoring an absorption identity gives additional factors, but does not by itself prove that these factors are nonisomorphic. Theorem~\ref{thm:map-nonunique-factors} below proves nonuniqueness for each fixed reflexive Lorentz pair by distinguishing factors according to the presence of $c_0$. Theorem~\ref{thm:bap-rigidity} concerns uniqueness only within the specified tensor orbit.} More generally, no regularity of the second factor is needed for collapse after $U$ has already been inserted.

\begin{corollary}\label{cor:principal-orbit}
Let $Z\ne\{0\}$ be any Banach space and put $X_Z=U\tp Z$.  For every pair satisfying \eqref{eq:banach-range},
\[
  \ell_{p,q}\tp X_Z\cong X_Z\cong L_{p,q}(0,1)\tp X_Z.
\]
\end{corollary}

\begin{proof}
By associativity and Theorem \ref{thm:universal-absorption},
\[
\begin{aligned}
 \ell_{p,q}\tp(U\tp Z)
 &\cong (\ell_{p,q}\tp U)\tp Z
 \cong U\tp Z,\\
 L_{p,q}(0,1)\tp(U\tp Z)
 &\cong (L_{p,q}(0,1)\tp U)\tp Z
 \cong U\tp Z.
\end{aligned}
\]
No approximation property, separability or basis hypothesis on $Z$ is used.
\end{proof}

\begin{remark}\label{rem:global-unique}
Theorem \ref{thm:bap-rigidity} is an orbit theorem.  It does not assert that every separable BAP factor $X$ satisfying
\[
  \ell_{p,q}\tp X\cong L_{p,q}(0,1)\tp X
\]
for all pairs \eqref{eq:banach-range} must be isomorphic to $U$.  Such a global uniqueness theorem would require a mechanism forcing an arbitrary repairing factor into the principal orbit of $U$, and no such mechanism is available here.
\end{remark}

\begingroup
\paragraph{Nonuniqueness among separable MAP factors.}
\begin{theorem}\label{thm:map-nonunique-factors}
Fix $1<p,q<\infty$, $p\ne q$, and put $A=\ell_{p,q}$ and $B=L_{p,q}(0,1)$. There exist nonzero separable Banach spaces $X,Y$, both having the metric approximation property (MAP)\footnote{A Banach space $E$ has the MAP if for every compact set $K\subset E$ and every $\varepsilon>0$, there is a finite-rank bounded linear operator $T:E\to E$ with $\|T\|\le 1$ and $\sup_{x\in K}\|Tx-x\|<\varepsilon$.}, such that
\[
 A\tp X\cong X\cong B\tp X,\qquad
 A\tp Y\cong Y\cong B\tp Y,\qquad X\not\cong Y.
\]
The space $X$ can be chosen to be a separable dual space containing no copy of $c_0$, whereas $Y$ contains a complemented copy of $c_0$.
\end{theorem}
\begin{proof}
Set $E_1=A$, $E_2=B$. These spaces are separable and reflexive. Equip them with the equivalent real interpolation $K$-norms for $(\ell_1,\ell_\infty)_{\theta,q}$ and $(L_1,L_\infty)_{\theta,q}$, respectively, with $\theta=1-1/p$. Coordinate truncations and conditional expectations onto finite dyadic partitions are finite-rank contractions for these norms. They converge strongly to the identity by density of finite sequences and dyadic step functions, respectively. A finite-net argument gives convergence uniformly on compact sets. Thus both generators have the MAP.

Let $\mathcal W$ consist of all finite words over $\{1,2\}$, including the empty word. Define
\[
 E_\varnothing=\K,\quad E_w=E_{i_1}\tp\cdots\tp E_{i_m}
 \quad(w=(i_1,\ldots,i_m)),\qquad
 X=\left(\bigoplus_{w\in\mathcal W,\ k\ge1}E_w\right)_{\ell_1}.
\]
This space is nonzero and separable. Every $E_w$ has the MAP: tensor products of finite-rank contractions are finite-rank contractions, and on elementary tensors their strong convergence follows from
\[
 \|(S\otimes T)(f\otimes g)-f\otimes g\|_\pi
 \le\|Sf-f\|\,\|g\|+\|f\|\,\|Tg-g\|.
\]
Density extends convergence to the completed tensor product. Finite-coordinate truncation followed by finite-rank contractions in the retained coordinates proves the MAP for $X$. All approximation constants are one, independently of word length.

Reindexing the infinite multiplicities gives $X\cong\ell_1(X)$. Distribution over $\ell_1$-sums identifies $E_i\tp X$ with the coordinate-complemented part of $X$ indexed by words beginning with $i$. Thus $E_i\tp X\ce X$. Choose norm-one $e_i\in E_i$, $e_i^*\in E_i^*$ with $e_i^*(e_i)=1$. The maps $J_i x=e_i\otimes x$ and $Q_i(e\otimes x)=e_i^*(e)x$ satisfy $Q_iJ_i=I_X$, giving $X\ce E_i\tp X$. The decomposition principle (equivalently Lemma~\ref{lem:r-decomp} with $r=1$) yields
\begin{equation}\label{eq:word-map-absorption}
 E_i\tp X\cong X\qquad(i=1,2).
\end{equation}

To distinguish the factors, we construct a predual of $X$. The classical tensor duality theorem states that
\[
 (F\widehat\otimes_\varepsilon G)^*=F^*\tp G^*
\]
isometrically if $F^*$ has the Radon--Nikod\'ym property and either $F^*$ or $G^*$ has the approximation property; see \cite[Theorem~5.33]{Ryan2002} and \cite[Theorem~5]{BotelhoFleming2021}. Put $D_i=E_i^*$, $G_\varnothing=\K$, and
\[
 G_w=D_{i_1}\widehat\otimes_\varepsilon\cdots\widehat\otimes_\varepsilon D_{i_m}.
\]
Induction gives $G_w^*=E_w$ isometrically. Indeed the empty-word case holds, and reflexivity gives $D_i^*=E_i$, which has the Radon--Nikod\'ym property and MAP. Hence the induction step is
\[
 G_{iw}^*=(D_i\widehat\otimes_\varepsilon G_w)^*
 =D_i^*\tp G_w^*=E_i\tp E_w=E_{iw}.
\]
The duality formula for $c_0$-sums now gives
\[
 X=\left(\left(\bigoplus_{w\in\mathcal W,\ k\ge1}G_w\right)_{c_0}\right)^*.
\]
Thus $X$ is a separable dual space. It has the Radon--Nikod\'ym property, which passes to closed subspaces and is invariant under isomorphism. Since $c_0$ lacks this property, $c_0\not\hookrightarrow X$.

Set $Y=X\tp c_0$. This is separable and has the MAP by the tensor approximation argument above. Associativity and \eqref{eq:word-map-absorption} give
\[
 E_i\tp Y\cong(E_i\tp X)\tp c_0\cong X\tp c_0=Y.
\]
For norm-one $x_0\in X$, $x_0^*\in X^*$ with $x_0^*(x_0)=1$, the maps $j(a)=x_0\otimes a$ and $r(x\otimes a)=x_0^*(x)a$ satisfy $rj=I_{c_0}$. Thus $c_0\ce Y$, and $X\cong Y$ would contradict $c_0\not\hookrightarrow X$.
\end{proof}
\begin{remark}\label{rem:nonunique-map-scope}
Together with \eqref{eq:base-nonembedding}, this gives two nonisomorphic separable MAP factors eliminating the original Lorentz nonembedding. The assertion concerns each fixed pair in the reflexive range; it does not address $q=1$ or uniqueness of a factor repairing every pair in \eqref{eq:banach-range} simultaneously. Moreover $X\not\cong U$, since the universal BAP absorber $U$ contains a complemented copy of $c_0$.
\end{remark}
\endgroup

\noindent\
The absorption results show that the Lorentz obstruction can disappear.  We now identify factors for which the opposite happens.  The proofs use finite self-similarity of Lorentz spaces and, for $q>1$, reflexivity of $\ell_{p,q}$; the finite-amplification and localization tools needed below were established in Section~\ref{sec:prelim}.

\begin{proposition}\label{prop:finite-barrier}
Assume \eqref{eq:banach-range}.  Let $Y\ne0$ be any Banach space and let $F\ne0$ be finite-dimensional.  Then
\[
    \ell_{p,q}\tp Y\not\hookrightarrow L_{p,q}(0,1)\tp F.
\]
In particular, every nonzero finite-dimensional common factor preserves \eqref{eq:base-nonembedding}.
\end{proposition}

\begin{proof}
If such an embedding existed, restrict it to the complemented rank-one copy of $\ell_{p,q}$ in $\ell_{p,q}\tp Y$.  Lemma \ref{lem:finite-amplification} identifies the target with $L_{p,q}(0,1)$ up to isomorphism, contradicting \eqref{eq:base-nonembedding}.
\end{proof}

\noindent\textbf{Stability classes and the $q=1$ threshold.}
The opposite phenomenon is that a large class of factors remains transparent to the Lorentz obstruction. We call $X$ an $\ell_1$-FDD factor if it is isomorphic to an
$\ell_1$-sum of finite-dimensional spaces, namely,
\begin{equation}
X\cong
\left(\bigoplus_{n=1}^{\infty}F_n\right)_{\ell_1},
\qquad
0<\dim F_n<\infty .
\end{equation}
This terminology combines the standard notion of a finite-dimensional
decomposition with the $\ell_1$-direct-sum norm; see
\cite[p.~198, Definition~4.31]{FabianHabalaHajekMontesinosZizler}
and \cite[Sec.~6, p.~298]{CasazzaApproximation}.
No uniform bound on $\dim F_n$ and no control of their Banach--Mazur
geometry is assumed.

\begin{theorem}\label{thm:l1fdd}
Assume $1<p<\infty$, $1<q<\infty$ and $p\ne q$.  If $X$ is a nonzero $\ell_1$-FDD factor, then for every nonzero Banach space $Y$,
\[
  \ell_{p,q}\tp Y\not\hookrightarrow L_{p,q}(0,1)\tp X.
\]
In particular,
\[
  \ell_{p,q}\tp X\not\hookrightarrow L_{p,q}(0,1)\tp X.
\]
\end{theorem}

\begin{proof}
Write $X\cong(\bigoplus_nF_n)_{\ell_1}$ with each $F_n$ finite-dimensional.  Suppose, toward a contradiction, that
\[
  S:\ell_{p,q}\tp Y\longrightarrow L_{p,q}(0,1)\tp X
\]
is an isomorphic embedding.  The rank-one retraction in the source gives an isomorphic copy of $\ell_{p,q}$, so by \eqref{eq:l1-distrib} we obtain
\begin{equation}\label{eq:into-l1sum}
 \ell_{p,q}\hookrightarrow
 \left(\bigoplus_{n=1}^{\infty}
 L_{p,q}(0,1)\tp F_n\right)_{\ell_1}.
\end{equation}
Because $1<p,q<\infty$, $\ell_{p,q}$ is reflexive.  Corollary \ref{cor:compression} compresses \eqref{eq:into-l1sum} into the first $N$ coordinates.  If
\[
        F=\left(\bigoplus_{n=1}^{N}F_n\right)_{\ell_1},
\]
then $F$ is finite-dimensional and \eqref{eq:l1-distrib} identifies the finite target with $L_{p,q}(0,1)\tp F$.  By Lemma \ref{lem:finite-amplification} this is isomorphic to $L_{p,q}(0,1)$.  We have therefore produced an embedding $\ell_{p,q}\hookrightarrow L_{p,q}(0,1)$, contradicting \eqref{eq:base-nonembedding}.
\end{proof}

\begin{corollary}\label{cor:separation}
Fix a pair \eqref{eq:banach-range} with $q>1$.  No nonzero $\ell_1$-FDD factor is isomorphic to a member of the principal absorbing orbit $\{U\tp Z:Z\ne0\}$.
\end{corollary}

\begin{proof}
Theorem \ref{thm:l1fdd} gives nonembedding for the first class, whereas Corollary \ref{cor:principal-orbit} gives isomorphism for the second.
\end{proof}

This separation is stronger than a dimension dichotomy.  Both classes contain infinite-dimensional spaces, and both are stable under many familiar operations.  What distinguishes them is the way their tensor structure localizes or absorbs complemented subspaces.

\begingroup
\begin{remark}\label{rem:q1threshold}
At $q=1$, the source $\ell_{p,1}$ is nonreflexive, so the unit-ball compression argument of Theorem~\ref{thm:l1fdd} is unavailable. The continuous inclusion $L_{p,1}(0,1)\subset L_1(0,1)$ does not supply a substitute isomorphic embedding. Indeed, for $f_t=t^{-1/p}\mathbf1_{(0,t)}$, $0<t<1$, one has
\[
 \|f_t\|_{L_{p,1}}\asymp1,\qquad
 \|f_t\|_{L_1}=t^{1-1/p}\longrightarrow0.
\]
Thus this inclusion is not bounded below, and cotype cannot be transferred through it as through an isomorphic embedding. The $\ell_1$-FDD preservation assertion at $q=1$ is not established here. Finite-dimensional preservation and universal absorption remain covered by Proposition~\ref{prop:finite-barrier} and Theorem~\ref{thm:universal-absorption}, respectively.
\end{remark}
\endgroup

\begingroup
The preceding discussion isolates the reflexive range in which coordinate
compression is available.  We now apply this mechanism to scalar $L_1$
factors.  The following lemma records the tail estimate underlying the
atomic case in a form suited to the classification argument.

\begin{lemma}\label{lem:pdf-atomic-tail}
Let
\[
       W=\left(\bigoplus_{n\ge1}W_n\right)_{\ell_1}
\]
and let $P_N$ denote the projection onto the first $N$ coordinates.  If
$K\subset W$ is weakly compact, then
\[
       \sup_{z\in K}\norm{(I-P_N)z}\longrightarrow0
       \qquad(N\to\infty).
\]
\end{lemma}

\begin{proof}
This is the special case $B_n=W_n$ of Lemma~\ref{lem:uniform-tails}, which
allows the coordinate spaces to vary.  It is precisely the finite-coordinate
compactness statement used in the atomic part of the classification argument.
\end{proof}

The preceding lemma supplies the localization step for purely atomic factors.
It remains to determine when a scalar $L_1$ factor contains enough nonatomic
structure to remove the Lorentz obstruction.  The next theorem gives this
classification by combining the atomic compression argument with the
embedding construction for $L_1(0,1)$.

\begin{theorem}\label{thm:scalar-l1-classification}
Let $1<q<p<2$, let $(\Omega,\Sigma,\mu)$ be a standard probability space, and put
$X=L_1(\Omega,\mu)$.  The following assertions are equivalent:
\begin{enumerate}[label=\textup{(\roman*)}]
\item
$\ell_{p,q}\tp X\hookrightarrow L_{p,q}(0,1)\tp X$;
\item
the nonatomic part of $\mu$ has positive measure;
\item
$X\cong L_1(0,1)$.
\end{enumerate}
In particular, among scalar $L_1$-spaces of standard probability spaces,
$L_1(0,1)$ is the unique isomorphism class of factors for which the
Lorentz nonembedding in \eqref{eq:base-nonembedding} disappears in this
parameter range.
\end{theorem}

\begin{proof}
Put $I=(0,1)$, $H=L_1(I)$, $A=\ell_{p,q}$, and $B=L_{p,q}(0,1)$.
We first prove \textup{(ii)}$\Rightarrow$\textup{(iii)}.  Decompose the standard
probability space, modulo null sets, into its nonatomic part and its atomic
part.  If the nonatomic part has positive measure, the classification of
standard nonatomic finite measure spaces gives
$L_1(\Omega_{\mathrm{na}})\cong H$.  The atomic part contributes either the
zero space, a finite-dimensional $\ell_1$-space, or $\ell_1$, after the atom
masses have been absorbed into the coordinates.  Thus
\[
       X\cong H\oplus_1 F,
\]
where $F$ is one of these three spaces.  Partition $I$ into countably many
intervals of positive length.  Averaging over this partition gives a bounded
projection on $H$ whose range is isomorphic to $\ell_1$; consequently
$H\cong K\oplus_1\ell_1$ for its kernel $K$.  Since
$\ell_1\oplus_1F\cong\ell_1$, we obtain $H\oplus_1F\cong H$, proving
\textup{(iii)}.

To prove \textup{(iii)}$\Rightarrow$\textup{(i)}, define, for
$a=(a_n)\in A$,
\[
       (Da)(t)=\sum_{n\ge1}a_n\mathbf 1_{[n-1,n)}(t),
       \qquad t>0.
\]
The map $D$ preserves the distribution functions of the sequence and its
step-function realization, hence gives an isomorphic embedding
$A\hookrightarrow L_{p,q}(0,\infty)$.  For $1<q<p<2$, the embedding theorem
of Carothers and Dilworth \cite[Corollary~2]{CarothersDilworth1991}, followed
by the stable-law embedding of Bretagnolle, Dacunha-Castelle and Krivine
\cite{BretagnolleDacunhaCastelleKrivine1966}, gives an isomorphic embedding
$J:A\to H$.  Thus there are constants $0<c\le C<\infty$ such that
\begin{equation}\label{eq:scalar-l1-J}
       c\norm{a}_A\le \norm{Ja}_H\le C\norm{a}_A
       \qquad(a\in A).
\end{equation}
For a simple $A$-valued function $f$ on $I$, set
\[
       (\mathcal Jf)(s,t)=\bigl(J(f(s))\bigr)(t),
       \qquad (s,t)\in I^2.
\]
By Fubini's theorem and \eqref{eq:scalar-l1-J},
\[
c\|f\|_{L_1(I;A)}
\leq
\|\mathcal{T}f\|_{L_1(I^2)}
=
\int_I \|J(f(s))\|_H\,ds
\leq
C\|f\|_{L_1(I;A)}.
\]
Density therefore extends $\mathcal J$ to an isomorphic embedding
$L_1(I;A)\hookrightarrow L_1(I^2)$.  Since $I^2$ and $I$ are isomorphic
nonatomic probability spaces, $L_1(I^2)\cong H$.  If $b\in B$ is a unit vector,
then $h\mapsto[s\mapsto h(s)b]$ is an isometric embedding
$H\hookrightarrow L_1(I;B)$.  Using the canonical Bochner identifications
\[
       A\tp H\cong L_1(I;A),
       \qquad B\tp H\cong L_1(I;B)
\]
for the projective tensor product \cite[Section~2.3]{Ryan2002}, and using
$X\cong H$, we obtain \textup{(i)}.

It remains to prove \textup{(i)}$\Rightarrow$\textup{(ii)}.  Suppose that the
nonatomic part is null.  Then $X$ is isomorphic to an $\ell_1$-sum of
one-dimensional spaces, so it is a nonzero $\ell_1$-FDD factor.  Theorem
\ref{thm:l1fdd}, applied with $Y=X$, rules out the embedding in
\textup{(i)}.  Explicitly, a hypothetical embedding restricts to a copy of
$A$ in
$L_1(\mu;B)\cong(\bigoplus_nB)_{\ell_1}$; the reflexivity of $A$ and
Lemma~\ref{lem:pdf-atomic-tail} compress this copy into finitely many
coordinates, and Lemma~\ref{lem:finite-amplification} then produces
$A\hookrightarrow B$, contradicting \eqref{eq:base-nonembedding}.
\end{proof}

For example, for $p=3/2$ and $q=4/3$, an $\ell_1$ factor preserves the
nonembedding, whereas the nonatomic factor $L_1(0,1)$ produces an embedding.
The same nonatomic factor works simultaneously for every pair
$1<q<p<2$.
\endgroup

\section{Quasi-Banach stabilization}\label{sec:quasi}

The Banach arguments above rely on a completed tensor norm, complemented rank-one slices, and a decomposition principle.  In the quasi-Banach setting these ingredients must be reformulated rather than formally reused.  We first set up the $\ell_r$-projective tensor product and the two structural tools that replace the corresponding Banach identities.

\noindent\textbf{The $\ell_r$-projective tensor product.}
An $r$-norm, $0<r\le1$, satisfies
\[
   \norm{x+y}^r\le\norm{x}^r+\norm{y}^r.
\]
A complete $r$-normed space is an $r$-Banach space.  By the Aoki--Rolewicz theorem \cite{Aoki1942,Rolewicz1957}; see also \cite{KPR1984}, every quasi-Banach space admits an equivalent $r$-norm for some $r>0$.  For a finite collection of spaces we may decrease the exponent and work with a common $r$.

For $r$-Banach spaces $E,F$ and $u\in E\otimes F$, set
\begin{equation}\label{eq:pir-intro}
 \pi_r(u)=\inf\left\{
 \left(\sum_{k=1}^n\norm{x_k}^r\norm{y_k}^r\right)^{1/r}:
 u=\sum_{k=1}^n x_k\otimes y_k
 \right\}.
\end{equation}
Ansorena--Bello \cite[Definitions 20--21 and Proposition 21]{AnsorenaBello2022} associate topological tensor products to symmetric sequence quasi-norms.  For the $\ell_r$ choice, set
\[
     N_r(E,F)=\{u\in E\otimes F:\pi_r(u)=0\},
\]
and denote by $E\tpr{r}F$ the completion of $(E\otimes F)/N_r(E,F)$ under the induced $r$-norm.  Thus $\tpr{r}$ is notation, not an assertion that $\pi_r$ is nondegenerate on the algebraic tensor product.  The construction is an $r$-Banach functor and bounded operators tensor to bounded operators \cite{AnsorenaBello2022}.

For clarity, we record the two structural identities used in our absorber argument.  They are elementary for the $\ell_r$ choice but are worth isolating because both replace classical Banach-tensor identities.

\begin{proposition}\label{prop:r-calculus}
Let $E$ and $(F_n)$ be $r$-Banach spaces.
\begin{enumerate}[label=\textup{(\roman*)}]
\item There is a canonical onto isometry
\begin{equation}\label{eq:r-distrib}
 E\tpr{r}\left(\bigoplus_nF_n\right)_{\ell_r}
 \cong
 \left(\bigoplus_nE\tpr{r}F_n\right)_{\ell_r}.
\end{equation}
\item For $r$-Banach spaces $E,F,G$, the canonical associator extends to an isometric isomorphism
\[
   (E\tpr{r}F)\tpr{r}G\cong E\tpr{r}(F\tpr{r}G).
\]
\end{enumerate}
\end{proposition}

\begin{proof}
For \textup{(i)}, on algebraic tensors send
\[
  \sum_k x_k\otimes(y_{k,n})_n
  \longmapsto
  \left(\sum_kx_k\otimes y_{k,n}\right)_n.
\]
For every finite representation,
\[
 \sum_n\pi_r\left(\sum_kx_k\otimes y_{k,n}\right)^r
 \le \sum_{n,k}\norm{x_k}^r\norm{y_{k,n}}^r
 =\sum_k\norm{x_k}^r\norm{(y_{k,n})_n}_{\ell_r(F_n)}^r.
\]
Thus the map is contractive.  Conversely, represent each coordinate of a finitely supported vector with $r$-cost arbitrarily close to its $\pi_r$ value and place those representations in disjoint coordinates.  This gives the reverse inequality.  Hausdorffization and density yield \eqref{eq:r-distrib}.

For \textup{(ii)}, on $E\otimes F\otimes G$ both bracketings induce the same three-fold projective $r$-seminorm
\[
 \inf\left\{
 \left(\sum_{k=1}^n\norm{x_k}^r\norm{y_k}^r\norm{z_k}^r\right)^{1/r}:
 w=\sum_{k=1}^n x_k\otimes y_k\otimes z_k
 \right\}.
\]
One inequality is obtained by expanding intermediate tensors; the reverse follows by approximating each intermediate $\pi_r$-representation before taking the outer infimum.  The two Hausdorff completions are therefore canonically isometric.
\end{proof}

Equivalent common $r$-norms lead to isomorphic tensor products: if each factor quasi-norm is changed by fixed two-sided equivalence constants, then the corresponding $\pi_r$-seminorms change by the product of those constants.  In particular, the isomorphic conclusions below do not depend on the particular Aoki--Rolewicz renorming chosen.

\noindent\textbf{Rank-one slices and decomposition without local convexity.}

The point at which quasi-Banach duality enters is extremely precise.

\begin{lemma}\label{lem:r-rankone}
Let $E,F$ be $r$-Banach spaces and suppose $F^*\ne\{0\}$.  Choose $y_0\in F$ and $\varphi\in F^*$ with $\varphi(y_0)=1$.  Then
\[
      J:E\longrightarrow E\tpr{r}F,\qquad Jx=x\otimes y_0,
\]
is an isomorphic embedding with complemented range.
\end{lemma}

\begin{proof}
Functoriality gives a bounded slice operator
\[
      Q=\mathrm{Id}_E\tpr{r}\varphi:E\tpr{r}F\longrightarrow E\tpr{r}\K\cong E.
\]
The canonical scalar identification is isometric: the $r$-triangle inequality gives the lower estimate for every scalar tensor representation, while the one-term representation gives the upper estimate.  Since $QJ=\mathrm{Id}_E$, the map $J$ is bounded below and $JQ$ is a bounded projection onto $J(E)$.
\end{proof}

The next lemma is the exact substitute for the decomposition step.  Its proof is purely topological and uses no Hahn--Banach separation.

\begin{lemma}\label{lem:r-decomp}
Let $0<r\le1$.  Suppose $X,Y$ are $r$-Banach spaces such that
\[
   X\ce Y,\qquad Y\ce X,\qquad
   X\cong\left(\bigoplus_{n=1}^{\infty}X\right)_{\ell_r}.
\]
Then $X\cong Y$.
\end{lemma}

\begin{proof}
Write $Y\cong X\oplus A$.  Because $Y\ce X$ and $A\ce Y$, transitivity gives $A\ce X$, so $X\cong A\oplus B$ for some $r$-Banach space $B$.  Using the assumed self-similarity and only coordinate regrouping,
\[
\begin{aligned}
 X\oplus\ell_r(A)
 &\cong \ell_r(X)\oplus\ell_r(A)\\
 &\cong \ell_r(A\oplus B)\oplus\ell_r(A)\\
 &\cong \ell_r(A)\oplus\ell_r(B)\oplus\ell_r(A)\\
 &\cong \ell_r(A)\oplus\ell_r(B)
 \cong \ell_r(X)
 \cong X.
\end{aligned}
\]
Since $\ell_r(A)\cong A\oplus\ell_r(A)$, it follows that
\[
 X\oplus A
 \cong (X\oplus\ell_r(A))\oplus A
 \cong X\oplus(\ell_r(A)\oplus A)
 \cong X\oplus\ell_r(A)
 \cong X.
\]
Consequently $Y\cong X\oplus A\cong X$.
\end{proof}

Let $J$ be finite or countable and $(E_j)_{j\in J}$ a family of $r$-Banach spaces.  Let $W(J)$ be the set of finite words over $J$, including the empty word.  For $w=(j_1,\dots,j_m)$ define
\[
   E_w=E_{j_1}\tpr{r}\cdots\tpr{r}E_{j_m},
   \qquad E_\varnothing=\K,
\]
where associativity is understood through Proposition \ref{prop:r-calculus}.  Give every word countably infinite multiplicity and put
\begin{equation}\label{eq:Ar-def}
 \mathfrak A_r(E_j:j\in J)
 =\left(
 \bigoplus_{w\in W(J)}\ \bigoplus_{m=1}^{\infty}E_w
 \right)_{\ell_r}.
\end{equation}

\begin{theorem}\label{thm:quasi-absorber-intro}
Let $(E_j)_{j\in J}$ be a finite or countable family of $r$-Banach spaces, where $0<r\le1$, and suppose $E_j^*\ne\{0\}$ for every $j$.  There exists an $r$-Banach space $\mathfrak A_r(E_j:j\in J)$ such that
\[
  E_j\tpr{r}\mathfrak A_r(E_i:i\in J)
  \cong \mathfrak A_r(E_i:i\in J)
  \qquad(j\in J).
\]
No approximation property is assumed.  If all $E_j$ are separable, the absorber is separable.
\end{theorem}

\begin{proof}
Fix $j\in J$ and abbreviate $\mathfrak A=\mathfrak A_r(E_i:i\in J)$.  By Proposition \ref{prop:r-calculus},
\[
 E_j\tpr{r}\mathfrak A
 \cong
 \left(
 \bigoplus_{w\in W(J)}\ \bigoplus_{m=1}^{\infty}E_{jw}
 \right)_{\ell_r}.
\]
The right-hand side is a coordinate-complemented subspace of $\mathfrak A$, namely the part indexed by words beginning with $j$.  Hence
\[
        E_j\tpr{r}\mathfrak A\ce\mathfrak A.
\]
Since $E_j^*\ne0$, Lemma \ref{lem:r-rankone}, together with symmetry of the tensor product, gives the reverse complemented embedding
\[
        \mathfrak A\ce E_j\tpr{r}\mathfrak A.
\]
By construction every word in \eqref{eq:Ar-def} already occurs infinitely many times, so a coordinate reindexing gives
\[
        \mathfrak A\cong\left(\bigoplus_{n=1}^{\infty}\mathfrak A\right)_{\ell_r}.
\]
Lemma \ref{lem:r-decomp} now yields $E_j\tpr{r}\mathfrak A\cong\mathfrak A$.  Countability of the summands gives separability when each $E_j$ is separable.
\end{proof}

Hunt's duality theory for Lorentz spaces \cite{Hunt1966} shows that the continuous dual of $L_{p,q}(0,1)$ is nontrivial for $p\ge1$ when $0<q\le1$.  Therefore the theorem applies to the atomic and nonatomic Lorentz models simultaneously.

For the Lorentz application, fix $1\le p<\infty$ and $0<q<1$, so that the parameters lie in the quasi-Banach range. We use the function-space quasi-norm introduced in \eqref{eq:lorentz-quasinorm}; the sequence model $\ell_{p,q}$ is equipped with the corresponding discrete rearrangement expression. Choose equivalent $r$-norms on both spaces for one common Aoki--Rolewicz exponent $0<r\le1$. We define
\begin{equation}\label{eq:lorentz-absorber-notation}
 \mathfrak A_{p,q}^{(r)}
 :=\mathfrak A_r\bigl(\ell_{p,q},L_{p,q}(0,1)\bigr),
\end{equation}
using the construction in \eqref{eq:Ar-def} with $J=\{1,2\}$, $E_1=\ell_{p,q}$ and $E_2=L_{p,q}(0,1)$. Thus $\mathfrak A_{p,q}^{(r)}$ is the $\ell_r$-sum, with countably infinite multiplicity, of all finite $r$-projective tensor words in these two spaces, including the empty word $\K$. The subscript $(p,q)$ records the generating pair, and the superscript $(r)$ records the exponent used for the tensor products and direct sum; it is not a tensor power and need not equal $p$ or $q$. This notation refers to the chosen construction, without asserting uniqueness among all common absorbing factors.
\begin{corollary}\label{cor:quasi-lorentz}
Let $1\le p<\infty$ and $0<q<1$.  After equivalent renorming by a common $r$-norm, there is a separable $r$-Banach space $\mathfrak A_{p,q}^{(r)}$ such that
\[
  \ell_{p,q}\tpr{r}\mathfrak A_{p,q}^{(r)}
  \cong \mathfrak A_{p,q}^{(r)}
  \cong L_{p,q}(0,1)\tpr{r}\mathfrak A_{p,q}^{(r)}.
\]
\end{corollary}

\begin{proof}
Choose a common Aoki--Rolewicz exponent $r$ for $\ell_{p,q}$ and $L_{p,q}(0,1)$.  The coordinate functional $a\mapsto a_1$ is nonzero and continuous on $\ell_{p,q}$.  Hunt's Lorentz duality theorem \cite{Hunt1966} implies $L_{p,q}(0,1)^*\ne\{0\}$ when $p>1$ and $0<q\le1$, while for $p=1$ and $0<q\le1$ the dual is $L_\infty$; in particular it is nonzero throughout $1\le p<\infty$, $0<q<1$.  Apply Theorem \ref{thm:quasi-absorber-intro} to the two-space family.  Both spaces are separable for $q<\infty$, so the resulting absorber is separable.
\end{proof}

This result should be compared carefully with Theorem \ref{thm:universal-absorption}.  The Banach factor $U$ is universal across all nonzero separable BAP spaces, whereas $\mathfrak A_{p,q}^{(r)}$ is constructed from the particular Lorentz pair.  On the other hand, the quasi-Banach theorem assumes no approximation property at all.  The price of abandoning the BAP is pair dependence and the minimal duality hypothesis needed for rank-one retractions.

\noindent\textbf{Finite-dimensional stability and the $p<1$ obstruction.}

Finite-dimensional factors remain transparent under $\tpr{r}$.  If $F$ is $m$-dimensional with basis $(f_j)$ and biorthogonal functionals $(f_j^*)$, the map
\[
     (x_1,\dots,x_m)\longmapsto\sum_{j=1}^m x_j\otimes f_j
\]
is an isomorphism from a finite $\ell_r$-sum of copies of $E$ onto $E\tpr{r}F$: the upper estimate follows directly from $\pi_r$, and the lower estimate follows from the bounded slice maps $\mathrm{Id}_E\tpr{r}f_j^*$.  Together with Lemma \ref{lem:finite-amplification}, this gives the following conditional stability theorem.

\begin{theorem}\label{thm:quasi-finite}
Fix $0<p,q<\infty$ and suppose the underlying isomorphic nonembedding
\begin{equation}\label{eq:quasi-base}
     \ell_{p,q}\not\hookrightarrow L_{p,q}(0,1)
\end{equation}
is known for the pair under consideration.  Choose a common Aoki--Rolewicz exponent $r$ for the Lorentz spaces.
\begin{enumerate}[label=\textup{(\roman*)}]
\item If $X$ is a nonzero $r$-Banach space with $X^*\ne0$ and $F\ne0$ is finite-dimensional, then
\[
   \ell_{p,q}\tpr{r}X\not\hookrightarrow L_{p,q}(0,1)\tpr{r}F.
\]
\item For arbitrary nonzero finite-dimensional $F$ and $G$,
\[
   \ell_{p,q}\tpr{r}F\not\hookrightarrow L_{p,q}(0,1)\tpr{r}G.
\]
\end{enumerate}
\end{theorem}

\begin{proof}
For \textup{(i)}, Lemma \ref{lem:r-rankone} embeds a complemented copy of $\ell_{p,q}$ into $\ell_{p,q}\tpr{r}X$.  Any embedding as in the forbidden conclusion would therefore give
\[
     \ell_{p,q}\hookrightarrow L_{p,q}(0,1)\tpr{r}F.
\]
Finite-dimensional tensor reduction and Lorentz finite self-similarity identify the target with $L_{p,q}(0,1)$, contradicting \eqref{eq:quasi-base}.  Part \textup{(ii)} follows because every nonzero finite-dimensional space has a nonzero continuous functional.
\end{proof}

The hypothesis $X^*\ne0$ in the one-sided statement cannot be suppressed by a Hahn--Banach argument: Hahn--Banach separation is precisely what may fail in a nonlocally convex space.  For the nonatomic Lorentz model this failure becomes absolute below $p=1$.

\begin{proposition}\label{prop:trivial-dual}
Let $0<p<1$ and $0<q<\infty$.  Then
\[
       L_{p,q}(0,1)^*=\{0\}.
\]
Consequently, $L_{p,q}(0,1)$ has no nonzero finite-dimensional complemented subspace and admits no complemented rank-one slice.
\end{proposition}

\begin{proof}
The dual identity is a classical consequence of Hunt's description of the duals of $L(p,q)$ spaces on nonatomic measure spaces \cite{Hunt1966}; see also Cwikel \cite{Cwikel1975} for the surrounding weak-$L_p$ duality picture.  If a nonzero finite-dimensional subspace $F$ were complemented by a bounded projection $P$, composing $P$ with a nonzero element of $F^*$ would produce a nonzero continuous functional on $L_{p,q}(0,1)$, a contradiction.
\end{proof}

\begin{corollary}\label{cor:architecture-limit}
For $0<p<1$ and $0<q<\infty$, no absorption proof for the nonatomic Lorentz factor can be based on a bounded retraction of the form
\[
  L_{p,q}(0,1)\tpr{r}X \longrightarrow X
\]
obtained by slicing against a nonzero functional on $L_{p,q}(0,1)$.  In particular, the mutual-complemented-embedding argument used in Theorem \ref{thm:quasi-absorber-intro} cannot cross the threshold $p=1$.
\end{corollary}

\begin{remark}\label{rem:not-no-absorber}
Corollary \ref{cor:architecture-limit} is an obstruction to a method, not a no-go theorem for all possible tensor stabilization.  It does not rule out a common collapse factor based on a different quasi-Banach tensor functor, nor does it rule out an argument that avoids complemented one-dimensional slices.  Any such extension would have to exploit structure absent from the Banach-style absorber mechanism.
\end{remark}

\begin{remark}[A concrete illustration of absorption]
The universal absorber used above is not merely abstract: starting from
Pe\l czy\'nski's complementably universal space $P$, it is explicitly realized as
\[
U=\left(\bigoplus_{n=1}^{\infty}P\right)_{\ell_1},
\]
as in \eqref{eq:U-def}.  For a concrete Lorentz pair, take $(p,q)=(3,2)$.  The
classical nonembedding gives
\[
\ell_{3,2}\not\hookrightarrow L_{3,2}(0,1),
\]
whereas Theorem \ref{thm:universal-absorption} yields
\[
\ell_{3,2}\tp U\cong U\cong L_{3,2}(0,1)\tp U.
\]
Thus the original atomic--nonatomic obstruction disappears after tensor
stabilization by the same explicitly constructed factor $U$.  This remark is
included only to illustrate the absorption mechanism already established above,
rather than to introduce an additional result.
\end{remark}

\section{Further tensor stability consequences and conclusion}\label{sec:further}

The preceding sections establish the main Banach and quasi-Banach stabilization mechanisms.  We close by recording several consequences that organize the same-factor problem and the complementary preservation mechanism without repeating the principal absorption and rigidity theorems.

\begin{corollary}
\label{cor:canonical-repair-factor}
For every pair of nonzero Banach spaces $A$ and $B$,
\[
 A\hookrightarrow B\tp A,
 \qquad
 B\hookrightarrow A\tp B .
\]
Consequently, every non-embedding obstruction can be repaired by adding the
opposite space as a projective tensor factor.
\end{corollary}

\begin{proof}
Apply Lemma \ref{lem:rank-one-tensor-retraction} first with $(A,X)=(A,B)$ and then with $(A,X)=(B,A)$.
\end{proof}

By Lemma \ref{lem:finite-amplification}, finite-dimensional tensor factors produce only finite amplifications of the Lorentz model.  Thus, in the classical Banach range, such factors do not create the universal absorption mechanism.  At the opposite extreme, Theorem \ref{thm:universal-absorption} already supplies nonzero common factors $X$ for which
\[
 \ell_{p,q}\tp X\cong L_{p,q}(0,1)\tp X,
\]
so the same-factor obstruction can genuinely disappear under stabilization.

\begin{proposition}\label{prop:finite-power-stability}
Let $E$ be a nonzero Banach space admitting the BAP.  For every $n\in\mathbb N$,
\[
 E^{\widehat{\otimes}_{\pi} n}\widehat{\otimes}_{\pi}U\cong U.
\]
\end{proposition}

\begin{proof}
The case $n=1$ is Theorem \ref{thm:universal-absorption}.  Repeated application of the same absorption identity and associativity of the projective tensor product gives the general case.
\end{proof}

\noindent\textbf{Preservation of non-embedding obstructions.}
While some tensor factors can eliminate non-embedding barriers, others preserve the original obstruction.  Let $B$ be a nonzero Banach space and define the associated $\pi$-target-small class by
\[
\mathfrak A_{\pi}(B)
=
\{Y\neq\{0\}: B\widehat{\otimes}_{\pi}Y\hookrightarrow B\}.
\]
No complementedness is required in this definition.

\begin{lemma}
\label{lem:elementary-nonembedding-transfer}
Let \(A\) and \(B\) be Banach spaces such that
\(A\not\hookrightarrow B\). Let \(X\) be a nonzero Banach space and assume
that
\[
B\widehat{\otimes}_{\pi}Y\hookrightarrow B
\]
for some nonzero Banach space \(Y\). Then
\[
A\widehat{\otimes}_{\pi}X
\not\hookrightarrow
B\widehat{\otimes}_{\pi}Y .
\]
In particular, whenever \(X\in\mathfrak A_{\pi}(B)\),
\[
A\widehat{\otimes}_{\pi}X
\not\hookrightarrow
B\widehat{\otimes}_{\pi}X .
\]
\end{lemma}

\begin{proof}
Assume to the contrary that
\(A\widehat{\otimes}_{\pi}X\hookrightarrow
B\widehat{\otimes}_{\pi}Y\).
By the definition of the target-small class,
\(B\widehat{\otimes}_{\pi}Y\hookrightarrow B\).
Hence
\[
A\widehat{\otimes}_{\pi}X\hookrightarrow B .
\]
Choose a norm-one functional \(x^{*}\in X^{*}\) and \(x_{0}\in X\) with
\(x^{*}(x_{0})=1\). The canonical rank-one embedding and the associated
slice map yield a complemented copy of \(A\) inside
\(A\widehat{\otimes}_{\pi}X\). Therefore \(A\hookrightarrow B\), which
contradicts the hypothesis.
\end{proof}

\begin{remark}
This criterion separates two opposite tensor phenomena. Universal absorber
factors may collapse Lorentz distinctions, whereas target-small factors
preserve the original non-embedding obstruction.
\end{remark}

\begin{theorem}
\label{thm:operator-localization}
Let $B$ be a Banach space and
$D=(\bigoplus_{n\geq1}F_n)_{\ell_1}$ with finite-dimensional $F_n$.
If $T:R\to B\widehat{\otimes}_{\pi}D$ is weakly compact, then
\[
\left\|T-(I_B\widehat{\otimes}P_N)T\right\|\longrightarrow0 .
\]
In particular, if $R$ is reflexive, every bounded operator satisfies this
localization property.
\end{theorem}

\begin{proof}
Using the projective tensor decomposition,
\[
B\widehat{\otimes}_{\pi}D
\cong
\left(\bigoplus_{n\geq1}
(B\widehat{\otimes}_{\pi}F_n)\right)_{\ell_1}.
\]
The image of the unit ball under a weakly compact operator is relatively

weakly compact. Applying Lemma~\ref{lem:uniform-tails} to this image, with $B_n=B\widehat{\otimes}_{\pi}F_n$, gives the required norm convergence.
\end{proof}

\begin{proposition}
\label{prop:stabilization-class-closure}
Let $\mathcal C$ denote the class of Banach spaces which are isomorphic
to complemented subspaces of $\ell_1$-sums of finite-dimensional spaces.
Then $\mathcal C$ is stable under finite projective tensor products.
\end{proposition}

\begin{proof}
If
$X\ce D_1$ and $Y\ce D_2$ are complemented
embeddings into $\ell_1$-FDD spaces, tensoring the corresponding
inclusion--projection pairs gives a complemented embedding
\[
X\widehat{\otimes}_{\pi}Y
\hookrightarrow
D_1\widehat{\otimes}_{\pi}D_2 .
\]
The latter space is again an $\ell_1$-sum of finite-dimensional blocks,
after reindexing the tensor-product blocks.
\end{proof}

\begin{theorem}
\label{thm:directional-vanishing}
Assume that $B$ has the finite-dimensional self-similarity property
$B^n\cong B$ for every $n\in\mathbb N$. Let $X$ belong to the
stabilization class above. If a reflexive space $R$ does not embed into
$B$, then no operator from $R\widehat{\otimes}_{\pi}Y$ into
$B\widehat{\otimes}_{\pi}X$ can be an isomorphic embedding for a
nonzero Banach space $Y$.
\end{theorem}

\begin{proof}
Choose an $\ell_1$-FDD space $D=(\bigoplus_nF_n)_{\ell_1}$ containing $X$ as a complemented subspace.  If
\[
S:R\widehat{\otimes}_{\pi}Y\longrightarrow B\widehat{\otimes}_{\pi}X
\]
were an isomorphic embedding, restrict $S$ to the complemented rank-one copy of $R$ in $R\widehat{\otimes}_{\pi}Y$ and compose with the canonical embedding of $B\widehat{\otimes}_{\pi}X$ into $B\widehat{\otimes}_{\pi}D$.  This gives an isomorphic embedding $T:R\to B\widehat{\otimes}_{\pi}D$.  Since $R$ is reflexive, $T$ is weakly compact, so Theorem \ref{thm:operator-localization} yields $N$ such that $(I_B\widehat{\otimes}P_N)T$ is still bounded below.  Its range lies in $B\widehat{\otimes}_{\pi}F$, where $F=\bigoplus_{n=1}^NF_n$ is finite-dimensional.  Lemma \ref{lem:finite-amplification} identifies this target with a finite direct sum $B^m$, hence with $B$ by the assumed self-similarity.  Thus $R\hookrightarrow B$, contradicting the hypothesis.
\end{proof}

The preceding results summarize the stability and absorption mechanisms governing Lorentz space nonembeddings under projective tensor products. The emphasis is on the preservation or disappearance of the original Lorentz obstruction, rather than on general tensor-product identities.  In the Banach category, complemented universality and localization separate absorbing factors from stable factor classes, while in the quasi-Banach category the same strategy survives after $r$-renorming precisely to the extent permitted by nontrivial duality.

\section*{Acknowledgements}
Thanks to all the members of the Functional Analysis Research team of the College of Mathematics and Statistics of Anqing Normal University for their discussion and correction of the difficulties and errors encountered in this paper.

\end{document}